\documentclass[10pt]{elsarticle}
\usepackage[a4paper,margin=2.5cm]{geometry}
\pdfoutput=1
\usepackage{lipsum}
\usepackage{amsmath}
\usepackage{amsthm}
\usepackage{amssymb}
\usepackage{graphicx}
\usepackage{esint}
\usepackage{titlesec}
\usepackage{amsfonts}
\usepackage{verbatim}
\usepackage{xcolor}
\usepackage[utf8]{inputenc}
\usepackage{amsfonts}
\usepackage{amsmath}
\usepackage{color}
\usepackage{footmisc}
\usepackage{bbold}
\usepackage{amsthm}
\usepackage{appendix}
\usepackage{graphicx}
\usepackage{halloweenmath} 
\usepackage{curve2e} 
\usepackage{hyperref}
\usepackage{bm}
\usepackage{float}
\usepackage{cancel}

\newtheorem{theorem}{Theorem}

\hypersetup{
    colorlinks=true,
    linkcolor=black,
    citecolor = black,
    filecolor=magenta,      
    urlcolor=cyan
    }

\makeatletter

\@ifdefinable\SuCmathpictvertex{} 
\@ifdefinable\@SuC@reserved@dimen{\newdimen\@SuC@reserved@dimen}

\newenvironment*{@SuC@math@picture}[8]{%

  \def\SuCmathpictvertex{\circle*{#6}}%
  \setlength\unitlength{\fontdimen 22 #5\tw@}%
  \setlength\@SuC@reserved@dimen{#7\unitlength}%
  \kern\@SuC@reserved@dimen
  \@HwM@d@pict@strut{#2}%
  \picture(#3,#1)(#4,-1)%
    \roundcap
    \roundjoin
    \linethickness{#8\@HwM@thickness@units@for #5}%
}{%
  \endpicture
  \kern\@SuC@reserved@dimen
}
\newcommand*\@SuC@general@pict[9]{%
  \begin{@SuC@math@picture}%
            {#2}{#3}
            {#4}{#5}
            #6
            {#7}
            {#8}
            {#9}
    #1%
  \end{@SuC@math@picture}%
}
\newcommand*\@SuC@math@version@shunt[7]{%
  \@HwM@choose@thicknesses{\@SuC@general@pict {#1}{#2}{#3}{#4}{#5}#7}%
      {{.4}{.2}{}}
      {{0.5}{.25}{0.75}}
}

\newcommand*\DeclareNewSuCMathPict[6]{%
  \newcommand*{#1}{%
    \@HwM@general@ordinary@symbol
      {\@SuC@math@version@shunt {#6}{#2}{#3}{#4}{#5}}%
  }%
}

\makeatother

\DeclareNewSuCMathPict{\pentagon}
            {3}{1}  
            {4}{-2} 
{
    \polygon(90:2)(162:2)(234:2)(306:2)(378:2)%
    \put (90:2){\SuCmathpictvertex}%
    \put(162:2){\SuCmathpictvertex}%
    \put(234:2){\SuCmathpictvertex}%
    \put(306:2){\SuCmathpictvertex}%
    \put(378:2){\SuCmathpictvertex}%
}

\newcommand{\per}[1]{\mathrm{\text{per}}\left({#1}\right)}

\newtheorem{Conjecture}{Conjecture \ignorespaces}
\newtheorem{Proposition}{Proposition \ignorespaces}
\newtheorem{corollary}{Corollary \ignorespaces}

\newtheorem*{counterexample*}{counterexample}
\newtheorem{lemma}{Lemma}

\newcommand{\Hn}{\mathcal{H}_{n}}

\DeclareMathOperator{\Per}{per}

\makeatletter
\def\ps@pprintTitle{%
 \let\@oddhead\@empty
 \let\@evenhead\@empty
 \let\@oddfoot\@empty
 \let\@evenfoot\@empty
}
\makeatother

\newcommand{\comm}[1]{\ignorespaces}

\begin{document}

\begin{frontmatter}

\title{A logical implication between two conjectures on matrix permanents}

\address[inst1]{Centre for Quantum Information and Communication,  École polytechnique de Bruxelles, CP 165/59,\\ Universit\'e libre de Bruxelles, 1050 Brussels, Belgium}

\address[inst2]{Physikalisches Institut, Albert-Ludwigs-Universit\"at Freiburg,
Hermann-Herder-Stra{\ss}e 3, D-79104 Freiburg, Germany}

\address[inst3]{EUCOR Centre for Quantum Science and Quantum Computing,

Albert-Ludwigs-Universit\"at Freiburg, Hermann-Herder-Stra{\ss}e 3, D-79104 Freiburg, Germany}

\address[inst4]{International Iberian Nanotechnology Laboratory (INL), Av. Mestre José Veiga, 4715-330 Braga, Portugal}

\author[inst1]{Léo Pioge}
\author[inst1]{Kamil K. Pietrasz}
\author[inst1,inst2,inst3,inst4]{Benoit Seron}
\author[inst4]{Leonardo Novo}
\author[inst1]{Nicolas J. Cerf}

\begin{abstract}
We prove a logical implication between two old conjectures stated by Bapat and Sunder about the permanent of positive semidefinite matrices. Although Drury has recently disproved both conjectures, this logical implication yields a non-trivial link between two seemingly unrelated conditions that a positive semidefinite matrix may fulfill. As a corollary, the classes of matrices that are known to obey the first conjecture are then immediately proven to obey the second one. Conversely, we uncover new counterexamples to the first conjecture by exhibiting a previously unknown type of counterexamples to the second conjecture.
Interestingly, such a relationship between these two mathematical conjectures appears from considerations on their quantum physics implications.
\end{abstract}

\end{frontmatter}
\vspace{-1pt}


\section{Introduction} 
Let $\Hn$ be the set of $n \times n$ positive semidefinite Hermitian matrices and $S_{n}$ be the symmetric group. The permanent of a $n \times n$ matrix $A=[a_{i,j}]$ is defined as 
\vspace {-0.3cm}
\begin{equation}
\per{A}=\sum_{\sigma\in S_{n}}\prod_{i=1}^{n}a_{i,\sigma(i)},
\end{equation}
where the summation is over the $n!$ permutations $\sigma$ of $n$ elements. It is well known that $\per{A}\ge 0$ provided $A$ is positive semidefinite \cite{Minc}. For a given $A\in\Hn$, let $F_{A}=[f_{i,j}]$ be the $n\times n$ matrix defined from $A$ as 
\begin{equation}
f_{i,j}=a_{i,j}\, \per{A(i,j)},
\end{equation}
where $A(i,j)$ is the $(n-1) \times (n-1)$ submatrix of $A$ obtained by deleting the $i$-th row and $j$-th column of $A$. Note that the Laplace expansion formula for the permanent implies that $\per A = \sum_i f_{i,j}, \forall j$, as well as $\per A=\sum_j f_{i,j}, \forall i$. Thus, it is clear that $\per{A}$ is an eigenvalue of $F_{A}$, which is associated with an eigenvector with all equal entries. It can also be shown that $F_A$ is a positive semidefinite matrix, so $F_A \in\Hn$. We note that, except for its eigenvalue $\per A$, its $n-1$ other (possibly degenerate) nonnegative eigenvalues are connected to $A$ in a non-trivial way. 

\section{Conjectures}
In 1985, Bapat and Sunder conjectured that the following inequality holds: 
\begin{Conjecture}[Bapat-Sunder \cite{bapat1985majorization}]
\label{conj:bapatSunder 1}
If $A =[a_{i,j}] \in \Hn$ and $B =[b_{i,j}] \in \Hn$, then 
\begin{equation}
\label{eq:conj1}
    \Per\,(A\circ B) \leq \Per\,(A) \prod_{i=1}^n b_{i,i},
\end{equation}
where $A\circ B$ denotes the Hadamard (or entry-wise) matrix product, namely $(A\circ B)=[a_{i,j} b_{i,j}]$.
\end{Conjecture}
This conjecture is the permanental counterpart to Oppenheim’s inequality \cite{oppenheim1930inequalities} on determinants, which states that $\det ({A\circ B}) \geq \det({A}) \prod_{i=1}^n b_{i,i}$. One year later, Bapat and Sunder proposed another conjecture involving the matrix $F_A$ that is deduced from $A$, namely:
\begin{Conjecture}[Bapat-Sunder \cite{bapat1986extremal}]
\label{conj:bapatSunder 2}
If $A =[a_{i,j}] \in \Hn$, then\footnote{In fact, Conjecture~\ref{conj:bapatSunder 2} 
 is often stated as the assertion that $\per{A}$ is the largest eigenvalue of matrix $F_A=[f_{i,j}] \in \Hn$, but this can equivalently be written as inequality \eqref{eq:conj2} since $F_{A}$ always admits $\per{A}$ as an eigenvalue.} 
\begin{equation}
 \label{eq:conj2}
 \lambda_\mathrm{max}(F_A) \leq \Per\,(A),
\end{equation}
where $\lambda_\mathrm{max}(F_A)$ denotes the largest eigenvalue of the matrix $F_A=[f_{i,j}] \in \Hn$, which is deduced from matrix $A$ following $f_{i,j}=a_{i,j}\, \Per\,(A(i,j))$.
\end{Conjecture}

Although both conjectures had remained open for about three decades, Drury recently exhibited a counterexample to Conjecture~\ref{conj:bapatSunder 1}~\cite{drury2016}, and, two years later, to Conjecture~\ref{conj:bapatSunder 2}~\cite{drury2018}. This invalidates both conjectures, of course, but does not make them less interesting. Indeed, characterizing the restricted classes of matrices that satisfy one or the other conjecture remains an intriguing -- and open -- problem, with applications in quantum physics, e.g., quantum interferometry \cite{shchesnovich2016universality,bosonbunching,pioge}. On the mathematics side, several works have contributed to this direction \cite{zhang1989notes,Zhang2013,pate}. 
Let us present some known results about the characterization of the matrices $A$ and $B$ that satisfy Conjecture \ref{conj:bapatSunder 1}.

\begin{Proposition}
\label{prop:0}
If $A \in \Hn$ is a rank-one matrix, then inequality \eqref{eq:conj1} in Conjecture~\ref{conj:bapatSunder 1} holds for all $B \in \Hn$.
\end{Proposition}
\noindent \textit{Proof:} If $A$ is a rank-one matrix in $\Hn$, then $A$ can be expressed as $A=\bm{x}^{*}\bm{x}$, where $\bm{x} \in \mathbb{C}^{n}$ and $\bm{x}^{*}=(\bm{\bar{x}})^{T}$ is the conjugate transpose of $\bm{x}$. Let $E$ be the $n \times n$ matrix whose entries are all equal to $1$, and define $\tilde{B}=[\tilde{b}_{i,j}]$ such that $\tilde{b}_{i,j}=b_{i,j}/\sqrt{b_{i,i}\, b_{j,j}}$. Since the matrix $\tilde{B} \in \Hn$ and all its diagonal entries are 1, it is a correlation matrix. As a consequence, we have, $\operatorname{per}(\tilde{B})=\sum_{\sigma\in S_{n}}\prod_{i=1}^{n}\tilde{b}_{i,\sigma(i)}\leq\sum_{\sigma\in S_{n}}\prod_{i=1}^{n}|\tilde{b}_{i,\sigma(i)}|\leq\sum_{\sigma\in S_{n}}\prod_{i=1}^{n} 1 =\per{E}$. From this, it follows that:
\begin{align}
\per{A \circ B} &= \per{(\bm{x}^{*}\bm{x}} \circ B) =\operatorname{per}\big([\overline{x_i} b_{ij} x_j]\big) \nonumber \\
&= \prod_{i=1}^n |x_i|^2 b_{i,i}~ \per{\tilde{B}} \leq \prod_{i=1}^n |x_i|^2 b_{i,i}~ \per{E} = \per{A} \prod_{i=1}^n b_{i,i}. \qed 
\end{align}

\begin{Proposition}[Zhang \cite{zhang1989notes}]
\label{prop:1}

 If $A \in \Hn$ has nonnegative entries ($a_{i,j}\ge 0, \forall i,j$), then inequality \eqref{eq:conj1} in Conjecture~\ref{conj:bapatSunder 1} holds for all $B \in \Hn$.
\end{Proposition}
\noindent \textit{Proof:} Since Zhang’s proof is short, we reproduce it here. If $B\in \Hn$, then $|b_{i,\sigma(i)} |\leq|b_{i,i}b_{\sigma(i),\sigma(i)}|^{1/2}$, $\forall \sigma \in S_{n}$, this gives  $\prod_{i=1}^n | b_{i,\sigma(i)} | \leq \prod_{i=1}^n b_{ii}$.  Furthermore, since $A$ is chosen such that $a_{i,j} \ge 0,
\forall i,j$, it follows that
\begin{equation}
 \per{A \circ B} =  \sum_{\sigma \in S_n} \prod_{i=1}^n a_{i,\sigma(i)} b_{i,\sigma(i)} 
\leq \sum_{\sigma \in S_n} \prod_{i=1}^n  a_{i,\sigma(i)} 
\left| b_{i,\sigma(i)} \right|
\leq \sum_{\sigma \in S_n} \prod_{i=1}^n a_{i,\sigma(i)} b_{i,i}
= \per{A} \prod_{i=1}^n b_{i,i}.\qed
\end{equation}

Thus, colloquially, we may say that Conjecture~\ref{conj:bapatSunder 1} holds provided $A \in \Hn$ is a rank-one matrix (Proposition~\ref{prop:0}) or a matrix with nonnegative entries (Proposition~\ref{prop:1}). It is natural, at this point, to question whether similar properties could be true when restricting matrix $B$ in Conjecture~\ref{conj:bapatSunder 1} instead of matrix $A$. The counterpart to Proposition~\ref{prop:0} was proven by Zhang, namely:

\begin{Proposition}[Zhang \cite{Zhang2013}]
\label{prop:2}
If $B \in \Hn$ is a rank-one matrix, then inequality \eqref{eq:conj1} in Conjecture~\ref{conj:bapatSunder 1} holds for all $A \in \Hn$ (actually, the inequality is saturated).
\end{Proposition}
\noindent \textit{Proof:} As the proof is succinct, we also provide it for completeness. If $B$ is a rank-one matrix in $\Hn$, then $B$ can be expressed as $B=\bm{x}^{*}\bm{x}$, where $\bm{x} \in \mathbb{C}^{n}$, then
\begin{equation}
\per{A \circ B} = \per{A\circ(\bm{x}^{*}\bm{x}} ) =\operatorname{per}\big([\overline{x_i} a_{ij} x_j]\big) = \prod_{i=1}^n |x_i|^2 \per{A}  = \prod_{i=1}^n  b_{i,i}~ \per{A}. \qed 
\end{equation}

The counterpart to Proposition~\ref{prop:1} appears to be a plausible proposition too, but Zhang's proof strategy does not easily extend, so we may only conjecture it at this point:

\begin{Conjecture}
\label{conj:3}
If $B =[b_{i,j}] \in \Hn$ has nonnegative entries ($b_{i,j}\ge 0, \forall i,j$) and $A =[a_{i,j}] \in \Hn$, then  
\begin{equation}
\label{eq:conj3}
    \Per\,(A\circ B) \leq \Per\,(A) \prod_{i=1}^n b_{i,i}.
\end{equation}
\end{Conjecture}


Since all counterexamples to Conjecture~\ref{conj:bapatSunder 1} known as of today \cite{drury2016,bosonbunching,drury2017real} do not involve matrices $B$ with nonnegative entries only, it is very tempting to assume that Conjecture~\ref{conj:3} holds as well. This would also nicely complete Propositions \ref{prop:0}-\ref{prop:2}. However, an unexpected consequence of the main theorem of this paper will be to disprove Conjecture~\ref{conj:3} (see Section \ref{sec:reverse}).


\section{Logically connecting Conjectures 1 and 2}
\label{sec:connecting}

Conjectures~\ref{conj:bapatSunder 1} and \ref{conj:bapatSunder 2} are linked to various similar conjectures involving the permanent of positive semidefinite matrices, such as the permanent-on-top conjecture \cite{POT}, Chollet's conjecture \cite{chollet}, Pate's conjecture \cite{pate}, or the Lieb's permanent dominance conjecture \cite{lieb1966}. The reader may consult the reviews of Wanless \cite{Wanless} or Zhang \cite{zhang2016} for a comprehensive perspective on these conjectures. We note that proving logical implications between conjectures is a standard tool in this area. This is the approach we follow in this paper: our main result is to unveil a logical implication between  Conjectures~\ref{conj:bapatSunder 1} and \ref{conj:bapatSunder 2}. By making a detour through a question motivated by quantum physics \cite{shchesnovich2016universality,bosonbunching,pioge}, we demonstrate that any counterexample to Conjecture~\ref{conj:bapatSunder 2} can be used to construct a counterexample that invalidates Conjecture~\ref{conj:bapatSunder 1}. Although both conjectures are now known to be false, this yields a valuable logical implication between specific properties of matrices in $\Hn$: namely, if Conjecture~\ref{conj:bapatSunder 1} holds for some class of matrices, then Conjecture~\ref{conj:bapatSunder 2} must hold for these matrices too. Conversely, if Conjecture~\ref{conj:bapatSunder 2} is disproved for a class of matrices, then Conjecture~\ref{conj:bapatSunder 1} is also disproved for those matrices. Our main theorem can be stated as follows:

\begin{theorem}
\label{th:central-theorem}  If $A \in \Hn$ is chosen such that inequality \ref{eq:conj1} in Conjecture~\ref{conj:bapatSunder 1} is verified for all $B \in \Hn$, then $A$ must also verify inequality \ref{eq:conj2} in
Conjecture~\ref{conj:bapatSunder 2}.
\end{theorem} 
 
Before proving the theorem, we consider a lemma. We define a $n$-dimensional complex vector $\bm{v}$, with $\lVert \bm{v}\rVert_2 = 1$, and a $n \times n$ correlation matrix 
$B(\epsilon)$ such that  $B(\epsilon)=M^{*}(\epsilon)M(\epsilon)$, where $M(\epsilon)$ is the $2\times n$ matrix defined as 
\begin{equation}
\label{eq:Mmat}
M(\epsilon) =\left(\begin{array}{cccccccc}
1/\alpha_{1} & 1/\alpha_{2} & \cdots & 1/\alpha_{n} \\ 
\epsilon v_{1}/\alpha_{1}&\epsilon v_{2}/\alpha_{2} & \cdots & \epsilon v_{n}/\alpha_{n}\\ 
\end{array}\right). 
\end{equation}
 Here, $\epsilon$ is a real positive variable, $v_{i}$ is the $i$-th component of vector $\bm{v}$, and $\alpha_{i} = \sqrt{1+\epsilon^{2} |v_{i}|^{2}}$. Thus, $B(\epsilon)$ is a positive semidefinite matrix whose diagonal elements are all equal to $1$, so it is a correlation matrix.  Note that $B(0)=E$. With these definitions, we are now ready to state our lemma:

\begin{lemma}\ignorespaces
    \label{lem:lemma1}
    Consider a matrix $A\in \Hn$ and a correlation matrix $ B(\epsilon)= M^*(\epsilon)M(\epsilon) \in \Hn$, with $M(\epsilon)$ defined as above. For $\epsilon>0$, the Taylor series of $\Per\,(A\circ B(\epsilon))$ around $\epsilon=0$ is given by
    \begin{equation}
       \Per\,(A\circ B(\epsilon))= \Per\,(A)+ \epsilon^2(\bm{v}^* F_{A} \bm{v}-\Per\,(A)) + O(\epsilon^4),
    \end{equation}
    where the matrix $F_A\in \Hn$ is deduced from matrix $A$ following $f_{i,j}=a_{i,j}\, \per{A(i,j)} $.
\end{lemma}

\noindent \textit{Proof:} In order to prove the lemma, we perform a Taylor expansion of all entries $b_{i,j}(\epsilon)$ of the matrix $B(\epsilon)$ in $\epsilon=0$, namely
\begin{align}
    b_{i,j}(\epsilon)=\frac{1}{\alpha_{i}\alpha_{j}}(1+\epsilon^2 v_{i}^{*}v_{j})=
    1+\epsilon^{2}\left(v_{i}^{*}v_{j}-
    \frac{|v_{i}|^{2}}{2}-\frac{|v_{j}|^{2}}{2}\right)+O(\epsilon^4).
\end{align}
To simplify the notation, we rewrite this equation as
\begin{equation}
    b_{i,j}(\epsilon)= 1 + \epsilon^{2}x_{i,j} +O(\epsilon^4),
\end{equation}
where we have introduced the $n\times n$ matrix $X=[x_{i,j}$] such that $x_{i,j}=v_{i}^{*}v_{j}-
    |v_{i}|^{2}/2-|v_{j}|^{2}/2$. 
This implies that 
\begin{equation}
\label{MincEq}
    \per{A\circ B(\epsilon)}=\per{A \circ (E+\epsilon^{2}X+O(\epsilon^4)) }= \nonumber 
    \per{A +\epsilon^{2}(A\circ X) }+O(\epsilon^4),
\end{equation}
Using a formula from Minc \cite{Minc} for the permanent of the sum of two matrices, we obtain,  
\begin{equation}
\label{eq:expansion_permanent}
    \per{A \circ B(\epsilon) }= \per{A}+ 
    \epsilon^{2} \sum_{i,j}^{n}(A\circ  X)_{i,j}\, \per{A(i,j)}+O(\epsilon^4),
\end{equation}
where $A(i,j)$ is the submatrix obtained by deleting the $i$-th row and $j$-th column of $A$. 
Given the definition of matrix $F_A=[f_{i,j}]$, with $f_{i,j}=a_{i,j}\, \per{A(i,j)}$, the coefficient of the term proportional to $\epsilon^2$ in Eq.~\eqref{eq:expansion_permanent} can be written as
\begin{equation}
     \sum_{i,j}^{n}x_{i,j}f_{i,j}= 
     \sum_{i,j}^{n} \Big(v_{i}^{*}v_{j}-
    \frac{|v_{i}|^{2}}{2}-\frac{|v_{j}|^{2}}{2}\Big)f_{i,j}=
     \sum_{i,j}^{n} v_{i}^{*}f_{i,j}v_{j}  - \sum_{i,j}^{n} \Big(\frac{|v_{i}|^{2}}{2}f_{i,j}-\frac{|v_{j}|^{2}}{2}f_{i,j}\Big).
\end{equation}
By using $\lVert \bm{v}\rVert_2 = 1$ and the Laplace expansion for the permanent, namely $\sum_{i}^{n} f_{i,j} = \sum_{j}^{n} f_{i,j} = \per{A}$, we have 
 \begin{equation}
     \sum_{i,j}^{n}x_{i,j}f_{i,j}= \bm{v}^{*}F_{A}\bm{v} - \per{A} .
 \end{equation}
By replacing this expression into Eq.~\eqref{eq:expansion_permanent}, we find that  
\begin{equation}
\label{eq:thm_statement}
    \per{A\circ B(\epsilon)}= \per{A}+ 
    \epsilon^2(\bm{v}^* F_{A} \bm{v}-\per{A}) + O(\epsilon^4), 
\end{equation}
which concludes the proof of the lemma. \qed

The implication relation between Conjectures~\ref{conj:bapatSunder 1} and \ref{conj:bapatSunder 2} that is contained in Theorem~\ref{th:central-theorem} can now be straightforwardly proven based on Lemma~\ref{lem:lemma1} as follows:


\medskip
\noindent \textbf{Proof of Theorem~\ref{th:central-theorem}:} 
Let $A\in \Hn$ be such that $ \lambda_\mathrm{max}(F_A)>\per{A}$, thus violating inequality \eqref{eq:conj2} in Conjecture \ref{conj:bapatSunder 2}. By choosing $\bm{v}$ to be the eigenvector corresponding to $ \lambda_\mathrm{max}(F_A)$,  we have $\bm{v}^* F_{A} \bm{v} > \per{A}$. Thus, for some $\epsilon>0 $, Lemma~\ref{lem:lemma1} implies that
\begin{equation}
    \per{A\circ B(\epsilon)} >\per{A}.
\end{equation}
Therefore, if matrix $A$ is a counterexample to Conjecture~\ref{conj:bapatSunder 2} [i.e., $A$ and $\bm{v}$ such that $\bm{v}^{*}F_{A}\bm{v}  > \per{A}$], then there exists some matrix $B(\epsilon)$ with $\epsilon\in (0, \epsilon_{max})$ such that $A$ and $B(\epsilon)$ yield a counterexample to Conjecture~\ref{conj:bapatSunder 1} [i.e., $\per{A\circ B(\epsilon)}>\per{A}$]. The contradiction of this logical implication concludes the proof of Theorem~\ref{th:central-theorem}. \qed

In short, Theorem~\ref{th:central-theorem} yields the direct implication Conjecture~\ref{conj:bapatSunder 1} $\Rightarrow$ ~\ref{conj:bapatSunder 2}, which was not known in the literature (see, e.g., \cite{Wanless,zhang2016}) and which will be used in Section~\ref{sec:direct}. Note that the reverse implication Conjecture~2 $\Rightarrow$ 1 is not valid, as can easily be tested with the example from Ref.~\cite{drury2016}, which obeys Conjecture~\ref{conj:bapatSunder 2} but contradicts Conjecture~\ref{conj:bapatSunder 1}. Instead, what will be exploited in Section~\ref{sec:reverse} is the contradictive implication ($\neg \,$Conjecture~\ref{conj:bapatSunder 2} $\Rightarrow \neg \,$Conjecture \ref{conj:bapatSunder 2}).

\section{Consequences of the direct logical implication}
\label{sec:direct}

Theorem~\ref{th:central-theorem} has nice ramifications as it provides some information on the eigenspectrum of matrix $F_A$ from existing knowledge on matrix $A$. In particular, any class of matrices $A$ known to satisfy Conjecture~\ref{conj:bapatSunder 1} (for all matrices $B$) is necessarily also a class of matrices $A$ satisfying Conjecture~\ref{conj:bapatSunder 2}. We have two interesting instances of this implication, which imply the following two corollaries of Theorem \ref{th:central-theorem}:

\begin{corollary}
\label{cor:1}
All rank-one positive semidefinite matrices $A$ satisfy inequality \eqref{eq:conj2} in Conjecture~\ref{conj:bapatSunder 2}.
\end{corollary}
By Proposition~\ref{prop:0}, it is known that all rank-one matrices $A$ satisfy Conjecture~\ref{conj:bapatSunder 1} for all $B \in \Hn$. Hence, we infer that rank-one matrices $A$ must satisfy Conjecture~\ref{conj:bapatSunder 2} too. This can also be verified directly as, if $A$ is of rank-one, then $F_A= \per{A} E/n$, hence its eigenvalues are 0 (with degeneracy $n-1$) and $\per{A}$.

\begin{corollary}
\label{cor:2}
All positive semidefinite matrices $A\in \Hn$ with nonnegative entries ($a_{i,j}\ge 0, \forall i,j$) satisfy inequality \eqref{eq:conj2} in Conjecture~\ref{conj:bapatSunder 2}.
\end{corollary}
By Proposition~\ref{prop:1} \cite{zhang1989notes}, it is known that all matrices $A \in \Hn$ with nonnegative entries satisfy Conjecture~\ref{conj:bapatSunder 1} for all $B \in \Hn$, hence they must satisfy Conjecture~\ref{conj:bapatSunder 2} too. This corollary allows us to easily recover a result that had been found earlier by Pate \cite{pate}.\\

Thus, Corollaries \ref{cor:1} and \ref{cor:2} illustrate the fact that the direct implication contained in Theorem \ref{th:central-theorem} yields non-trivial results about classes of matrices that satisfy Conjecture~\ref{conj:bapatSunder 2}. Although Corollaries \ref{cor:1} and \ref{cor:2} confirm known results, we expect that new results could be obtained following the same procedure. Furthermore, as we shall see in the next Section, it so happens that the contradictive implication brought by Theorem \ref{th:central-theorem} yields previously unknown (classes of) counterexamples to Conjecture~\ref{conj:bapatSunder 1}.


\section{Consequences of the contradictive logical implication} 
\label{sec:reverse}


Here, we will disprove Conjecture~\ref{conj:3} by introducing still another conjecture (Conjecture \ref{conj:4}) that is implied by Conjecture~\ref{conj:3} (in analogy with Conjecture~\ref{conj:bapatSunder 1} $\Rightarrow$ \ref{conj:bapatSunder 2}) and then by exploiting the contradictive implication. As we shall see, the validity of Conjecture~\ref{conj:bapatSunder 1} in the restricted case of $B$ matrices having nonnegative entries (i.e., what we called Conjecture~\ref{conj:3}) implies the following conjecture: 
\begin{Conjecture}
\label{conj:4}
If $A =[a_{i,j}] \in \Hn$,  then 
\begin{equation}
\label{eq:conj4}
\lambda_{\max}^{\mathbb{R}}(F_{A}) \leq \Per\,(A),
\end{equation}
where $F_A=[f_{i,j}]$ is deduced from matrix $A$ following $f_{i,j}=a_{i,j}\, \per{A(i,j)}$, and $\lambda_{\max}^{\mathbb{R}}(F_{A})$ stands for the maximum eigenvalue of the symmetric part of $F_A$ over $\mathbb{R}^{n\times n}$, i.e., the matrix formed by taking the real part of each entry of $F_A$.
\end{Conjecture}

While Conjecture~\ref{conj:3} is nothing but a restricted case of Conjecture~\ref{conj:bapatSunder 1},
Conjecture~\ref{conj:4} is close to -- but yet distinct from -- Conjecture~\ref{conj:bapatSunder 2} as it refers to the symmetric part of $F_A$ over $\mathbb{R}^{n\times n}$ instead of $F_A$. For this reason, we need to slightly adapt Theorem~\ref{th:central-theorem}, leading to the following theorem: 
\begin{theorem}
\label{th:2}
     If $A \in \Hn$ is chosen such that inequality \ref{eq:conj3} in Conjecture~\ref{conj:3} is verified for all $B \in \Hn$ with nonnegative entries ($b_{i,j}\ge 0, \forall i,j$), then $A$ must also verify inequality \ref{eq:conj4} in  Conjecture~\ref{conj:4}.
\end{theorem}

\medskip
\noindent \textbf{Proof of Theorem 2:} By performing a Taylor expansion of all entries $b_{i,j}(\epsilon)$ of the matrix $B(\epsilon)$ at $\epsilon=0$, we have $b_{i,j}(\epsilon)=
    1+\epsilon^{2}\left(v_{i}^{*}v_{j}-
    \frac{|v_{i}|^{2}}{2}-\frac{|v_{j}|^{2}}{2}\right)+O(\epsilon^4)$. If $\epsilon$ is sufficiently close to zero, and $v_i\in \mathbb{R}, ~\forall i$, then $B(\epsilon)$ has nonnegative entries. If $A$ and $B(\epsilon)$ fulfill Conjecture \ref{conj:bapatSunder 1}, then by Lemma~\ref{lem:lemma1} we have $\bm{v}^* F_{A} \bm{v}\leq\per{A}$, for any $\bm{v}\in \mathbb{R}^n$. Since $F_{A}\in \Hn $, this last condition can be rewritten as $ \lambda_{\max}^{\mathbb{R}}(F_{A})  \leq \per{A}$ since the maximum eigenvalue of the symmetric part of $F_A$ over $\mathbb{R}^{n\times n}$ is defined as the Rayleigh quotient over $\mathbb{R}^n$, that is, 
$\lambda_{\max}^{\mathbb{R}}(F_{A}) = \max_{\substack{\bm{v}\in\mathbb{R}^n\\ \lVert \bm{v}\rVert_2 = 1}} \bm{v}^T F_A \bm{v}$.\qed
\medskip 

In short, Theorem \ref{th:2} means that Conjecture~\ref{conj:3} implies Conjecture \ref{conj:4}, in close analogy with Theorem~\ref{th:central-theorem}. As a corollary of Theorem \ref{th:2}, we now exploit its reverse and thereby disprove Conjecture~\ref{conj:3}. This is our last result:

\begin{corollary}
\label{cor:3}
     Conjecture~\ref{conj:3} can be falsified, that is, it is not sufficient to consider a matrix $B\in\Hn$ with nonnegative entries ($b_{i,j}\ge 0, \forall i,j$) to ensure that inequality \eqref{eq:conj1} in Conjecture~\ref{conj:bapatSunder 1} holds for all $A \in \Hn$.
\end{corollary}

\medskip
\noindent \textit{Proof:} We start from a numerically found counterexample to Conjecture~\ref{conj:4}, which is a matrix $A$ of size
$16\times16$ constructed as $A=X^{*}X$, with $X=X_{R}+iX_{I}$ and
\begin{align}
X_{R}&=\left(\begin{array}{ccccccccccccccccc}
 25 & -23 &  29 &  11 & -20 & 47 &  18 & 29 & 35 & -25 & -32 & -28 & -18 & 25 & 12 & -36 \\
  8 &  38 & -11 &  34 &  61 & 42 & -23 & 10 & 35 &  24 &  11 &   9 &  13 &   -9 & 34 &  22
\end{array}\right) \nonumber\\
X_I&=\left(\begin{array}{ccccccccccccccccc}
 30 &  20 &  51 & -43 & -11 & 47 &   4 & 27 &  -26 &  -2 &  11 &  37 &  64 & 26 &  -28 &  23 \\
   0 &  20 &  10 &   4 &  28 & 12 & -46 & 24 &  -43 &  10 & -17 & -63 & -23 & 50 &  -40 &  15
\end{array}\right). 
\label{eq:XR_XI}
\end{align}
A numerical evaluation gives
\begin{align}
    &\lambda_{\max}^{\mathbb{R}}(F_{A})  \approx 2.2632\times10^{64}, \nonumber \\
    &\per{A}\approx 2.1978 \times10^{64},
\end{align}
which yields a ratio $\lambda_{\max}^{\mathbb{R}}(F_{A})/ \per{A}$ of approximately $1.0298$\footnote{This ratio may be lifted up to 1.0300 by optimizing the choice of $A$ in the neighborhood of the instance of $X_R$ and $X_I$ corresponding to Eq. \eqref{eq:XR_XI}, but we have chosen to build $X_R$ and $X_I$ with integer numbers for visual clarity.}, thereby providing a clear violation to Conjecture~\ref{conj:4}. Consequently, by the contrapositive of Theorem~\ref{th:2},
there exist some matrices $B(\epsilon)$ with nonnegative entries, where $\epsilon\in (0, \epsilon_{max})$ such that $A$ and $B(\epsilon)$ yield a counterexample to Conjecture~\ref{conj:3}, where $B(\epsilon)=M^{T}(\epsilon)M(\epsilon)$ is constructed according to 
Eq.~(\ref{eq:Mmat}), using as vector $\bm{v}$ the eigenvector of the symmetric part of $F_A$ over $\mathbb{R}^n$, associated to $\lambda_{\max}^{\mathbb{R}}(F_{A})$.
\qed 

Note that thorough numerical investigations have been performed for $n = 15$ but no counterexample was found, hence the smallest currently known counterexamples to Conjectures~\ref{conj:3} and \ref{conj:4} are sixteen-dimensional. Let us also mention that Theorem~\ref{th:2} was crucial in this numerical search for counterexamples to Conjecture~\ref{conj:3} with large matrices because finding counterexamples to Conjecture~\ref{conj:4} involves less parameters to optimize over (Conjecture~\ref{conj:4} involves a single matrix $A$ whereas Conjecture~\ref{conj:3} involves two matrices $A$ and $B$).

\section*{Acknowledgements}
We would like to thank Ian Wanless for correspondence and valuable suggestions.  L.P. is a FRIA grantee of the Fonds de la Recherche Scientifique – FNRS. B.S. was a research fellow of Fonds de la Recherche Scientifique  – FNRS and acknowledges funding from the Georg H. Endress Foundation. L.N. acknowledges funding from FCT-Fundação para a Ciência e a Tecnologia (Portugal) via the Project No. CEECINST/00062/2018. N.J.C. acknowledges funding from the Fonds de la Recherche Scientifique – FNRS under project
CHEQS (Grant No. 40007526) within the Excellence of Science (EOS) program.
L.N. and B.S. acknowledge support from Horizon Europe project EPIQUE (Grant No. 101135288).
\section*{Disclosure statement}
No potential conflict of interest was reported by the author(s).

\begin{center}

\bibliographystyle{unsrt}
\bibliography{main}

\end{center}

\end{document}